\documentclass[12pt]{article}
\usepackage{amsmath,amsthm,amsfonts,amssymb}
\usepackage[pdftex,pdfborder={0 0 0}]{hyperref}
\usepackage{fullpage}
\usepackage{multirow}
\usepackage{array}
\usepackage{bbm}
\usepackage[noblocks]{authblk}
\usepackage{verbatim}
\usepackage{tikz}
\usepackage{float}
\usepackage{enumerate}
\usepackage{diagbox}
\usepackage{comment}
\usetikzlibrary{arrows}
\usepackage{rotating}
\usepackage{soul}

\newtheorem{theorem}{Theorem}[section]

\newtheorem{conjecture}[theorem]{Conjecture}

\theoremstyle{definition}

\newlength{\Oldarrayrulewidth}

\newcommand{\Z}{\mathbb{Z}}

\newcommand{\p}{\textup{\textbf{p}}}
\newcommand{\q}{\textup{\textbf{q}}}

\newcommand{\lcm}{\textup{lcm}}

\renewcommand{\mod}[2]{\equiv#1\textup{ (mod }#2\textup{)}}

\def\m@th{\mathsurround=0pt}
\def\sm#1{\null\,\vcenter{\baselineskip9pt\lineskip.23ex\m@th
    \ialign{\hfil$\scriptstyle##$\hfil&&\ \hfil$\scriptstyle##$\hfil\crcr
    \mathstrut\crcr\noalign{\kern-\baselineskip}
    #1\crcr\mathstrut\crcr\noalign{\kern-\baselineskip}}}\,}
\def\smnp#1{\null\,\vcenter{\baselineskip9pt\lineskip.23ex\m@th
    \ialign{\hfil$\scriptstyle##$\hfil&&\ \ \hfil$\scriptstyle##$\hfil\crcr
    \mathstrut\crcr\noalign{\kern-\baselineskip}
    #1\crcr\mathstrut\crcr\noalign{\kern-\baselineskip}}}\,}

\begin{document}

\title{Probabilistic chip-collecting games with modulo winning conditions}

\author[1]{Joshua Harrington\thanks{joshua.harrington@cedarcrest.edu}}
\affil[1]{Department of Mathematics, Cedar Crest College}

\author[2]{Xuwen~Hua\thanks{xhaa2019@mymail.pomona.edu}}
\affil[2]{Department of Mathematics, Pomona College}

\author[3]{Xufei~Liu\thanks{xliu725@gatech.edu}}
\affil[3]{Department of Industrial Engineering, Georgia Institute of Technology}

\author[4]{Alex~Nash\thanks{nasha@dickinson.edu}}
\affil[4]{Department of Mathematics, Dickinson College}

\author[5]{Rodrigo~Rios\thanks{rodrigoreyrios@gmail.com}}
\affil[5]{Department of Mathematics, Florida Atlantic University}

\author[6]{Tony~W.~H.~Wong\thanks{wong@kutztown.edu}}
\affil[6]{Department of Mathematics, Kutztown University of Pennsylvania}

\date{\today}

\maketitle

\begin{abstract}

Let $a$, $b$, and $n$ be integers with $0<a<b<n$.  In a certain two-player probabilistic chip-collecting game, Alice tosses a coin to determine whether she collects $a$ chips or $b$ chips.  If Alice collects $a$ chips, then Bob collects $b$ chips, and vice versa.  A player is announced the winner when they have accumulated a number of chips that is a multiple of $n$.  In this paper, we settle two conjectures from the literature related to this game. 
\end{abstract}

\section{Introduction}\label{sec:intro}

In a probabilistic chip-collecting game introduced by Wong and Xu \cite{wx}, Alice and Bob take turns to toss a coin with Alice tossing first, which determines independently whether the player collects $a$ chips or $b$ chips.  The winner of the game is the first player to accumulate $n$ chips.  Some variations of this game have been considered by Leung and Thanatipanonda~\cite{lt,lt2020} and Harrington et al.~\cite{hkksw}.  The versions of the game that were considered by Harrington et al.\ removed the independence of the chip collecting process, so that if Alice collects $a$ chips, then Bob collects $b$ chips, and vice versa.  In one of these versions, called the \emph{modulo dependent game}, a player is announced the winner when they have accumulated a number of chips that is a multiple of $n$.

For $a<b<n$, the modulo dependent game can be treated as a random walk on $\mathbb{Z}_n\times\mathbb{Z}_n$, where the number of chips accumulated by each player is recorded as an ordered pair $(x,y)$ and each move is represented by either $(+a,+b)$ or $(+b,+a)$.  Since Alice always collects chips first, for any $y\in\Z_n$ and $x\in\Z_n\setminus\{0\}$, positions $(0,y)$ and $(x,0)$ are called the \emph{winning positions} of Alice and Bob, respectively, and a random walk on $\Z_n\times\Z_n$ that starts from $(0,0)$ \emph{terminates} upon landing on any winning position.  A position $(x,y)\in\mathbb{Z}_n\times\mathbb{Z}_n$ is said to be \emph{reachable} if there exists a random walk that lands on $(x,y)$ after leaving the starting position $(0,0)$.  As established by Harrington~et~al.~\cite{hkksw}, $(a,a)$ and $(b,b)$ are never reachable in $\Z_n\times\Z_n$.  They further conjectured the following statement, for which we provide a proof in Section~\ref{sec:section2}.

\begin{theorem}\label{thm:theorem1}
Every position in $\mathbb{Z}_n\times\mathbb{Z}_n\setminus\{(a,a),(b,b)\}$ is reachable if and only if $a\not\equiv 2b\pmod{n}$, $2a\not\equiv b\pmod{n}$, and $b^2-a^2$ is relatively prime to $n$.
\end{theorem}

The modulo dependent game can naturally be extended to a variation that allows Alice and Bob to having different winning conditions.  In particular, Harrington et al.\ considered a variation of the game where Alice wins by collecting a multiple of $m$ chips and Bob wins by collecting a multiple of $n$ chips.  This game can be recognized as a random walk on $\mathbb{Z}_m\times\mathbb{Z}_n$, where $a<b<\min\{m,n\}$.  Although this variation was not studied by Harrington et al., they did present the following conjecture.

\begin{conjecture}\label{conj:conjecture2}
Let $m\mid n$.  If all winning positions are of the form $(0,y)$, then $m\mid(b^2-a^2)$.
\end{conjecture}

In Section~\ref{sec:section3}, we will prove the following theorem, which establishes Conjecture~\ref{conj:conjecture2}.

\begin{theorem}\label{thm:theorem2}
In the modulo dependent game with parameters $a$, $b$, $m$, and $n$ such that $\gcd(a,b,m,n)=1$, all reachable winning positions are of the form $(0,y)$ if and only if $m\mid(b^2-a^2)$ and $m\mid\gcd(a,b)\gcd(m,n)$.
\end{theorem}

As a corollary to Theorem~\ref{thm:theorem2}, in the modulo dependent with parameters $a$, $b$, $m$, and $n$, notice that Bob's winning probability is $0$ if and only if $m\mid(b^2-a^2)$ and $m\mid\gcd(a,b)\gcd(m,n)$.


\section{Proof of Theorem~\ref{thm:theorem1}}\label{sec:section2}

\begin{proof}
If every position in $\Z_n\times\Z_n\setminus\{(a,a),(b,b)\}$ is reachable, then $(1,0)$ is reachable. In other words, $(ai+bj,aj+bi)=(1,0)$ for some integers $i$ and $j$. By adding or subtracting the two coordinates, we have $(a+b)(i+j)\equiv(b-a)(j-i)\mod{1}{n}$, thus $\gcd(a+b,n)=\gcd(b-a,n)=1$. Hence, $b^2-a^2$ is relatively prime to $n$. To establish the remaining necessary conditions, we proceed with a proof by contrapositive. If $a\mod{2b}{n}$, then the position $(2b,3b)$ can only be reached from $(0,2b)$ or $(b,b)$, so $(2b,3b)$ is not reachable. Similarly, if $b\mod{2a}{n}$, then the position $(2a,3a)$ is not reachable.

To prove the sufficient condition, let $\q_{i,j}=(ia+j(a+b),ib+j(a+b))$, where $i,j\in\Z$. Since $\gcd(b-a,n)=\gcd(a+b,n)=1$, every position in $\Z_n\times\Z_n$ can be expressed in the form of $\q_{i,j}$ for some $0\leq i,j\leq n-1$. Furthermore, $\gcd(k(a+b),n)\leq k<n$ and $\gcd(k(b-a),n)\leq k<n$ for all $1\leq k<n$, thus
\begin{equation}\label{eqn:k(a+b)}
k(a+b)\not\mod{0}{n}\text{ and }k(b-a)\not\mod{0}{n}.
\end{equation}
As a result, $2a\not\mod{2b}{n}$, which implies that every position $(x,x)\in\Z_n\times\Z_n\setminus\{(a,a),(b,b)\}$ is reachable by Harrington~et~al.~\cite[Theorem~3.6]{hkksw}. 
Hence, it remains to show that $\q_{i,j}$ is reachable for all $1\leq i\leq n-1$ and $0\leq j\leq n-1$.

We will prove by induction on $j$ that $\q_{1,j}$ is reachable for all $0\leq j\leq n-1$. First, the position $\q_{1,0}=(a,b)$ is reachable, and the position $\q_{1,1}$ is reachable by the sequence of moves
$$\q_{1,0}=(a,b)\xrightarrow{(+b,+a)}(a+b,a+b)\xrightarrow{(+a,+b)}\q_{1,1}.$$
Now, assume that for some $1\leq j\leq n-2$, $\q_{1,j'}$ is reachable for all $0\leq j'\leq j$. We proceed by considering the following cases. 
\begin{itemize}
\setlength{\itemindent}{15pt}
\item[\textit{Case $1$}:] $\q_{1,j}$ is not a winning position.

\begin{itemize}
\setlength{\itemindent}{30pt}
\item[\textit{Case $1(a)$}:] $\q_{2,j}$ is not a winning position.\\
The position $\q_{1,j+1}$ is reachable by the sequence of moves
$$\q_{1,j}\xrightarrow{(+a,+b)}\q_{2,j}\xrightarrow{(+b,+a)}\q_{1,j+1}.$$
\item[\textit{Case $1(b)$}:] $\q_{2,j}$ is a winning position.\\
Since $\q_{2,j}=(2a+j(a+b),2b+j(a+b))$, with a simple calculation, we have $\q_{2,j}\in\{(0,2b-2a),(2a-2b,0)\}$. Hence, $\q_{0,j+1}\in\{(b-a,b-a),(a-b,a-b)\}$, which does not intersect with $\{(a,a),(b,b)\}$ since $a\not\mod{2b}{n}$ and $b\not\mod{2a}{n}$. Therefore, $\q_{1,j+1}$ is reachable by the sequence of moves
$$\q_{1,j}\xrightarrow{(+b,+a)}\q_{0,j+1}\xrightarrow{(+a,+b)}\q_{1,j+1}.$$
\end{itemize}

\item[\textit{Case $2$}:] $\q_{1,j}$ is a winning position.\\
Since $\q_{1,j}=(a+j(a+b),b+j(a+b))$, with a simple calculation, we have $\q_{1,j}\in\{(0,b-a),(a-b,0)\}$. Hence, $\q_{1,j-1}\in\{(-a-b,-2a),(-2b,-a-b)\}$.

\begin{itemize}
\setlength{\itemindent}{30pt}
\item[\textit{Case $2(a)$}:] $\q_{1,j-1}$ is not a winning position.\\
Note that $b-2a\not\mod{0}{n}$ and $a-2b\not\mod{0}{n}$ by the given conditions, and $2b-2a\not\mod{0}{n}$ by \eqref{eqn:k(a+b)}. Hence, $\q_{2,j-1}\in\{(-b,b-2a),(a-2b,-a)\}$, $\q_{3,j-1}\in\{(a-b,2b-2a),(2a-2b,b-a)\}$, and $\q_{2,j}\in\{(a,2b-a),(2a-b,b)\}$ are not winning positions. Therefore, $\q_{1,j+1}$ is reachable by the sequence of moves
$$\q_{1,j-1}\xrightarrow{(+a,+b)}\q_{2,j-1}\xrightarrow{(+a,+b)}\q_{3,j-1}\xrightarrow{(+b,+a)}\q_{2,j}\xrightarrow{(+b,+a)}\q_{1,j+1}.$$
\item[\textit{Case $2(b)$}:] $2a\mod{0}{n}$ and $\q_{1,j-1}=(-a-b,0)$.\\
Note that $j>1$ since $\q_{1,1-1}=(a,b)\neq(-a-b,0)$. Also note that $-2a-2b\equiv-2b\not\equiv-2a\mod{0}{n}$ and $-a-2b\equiv a-2b\not\mod{0}{n}$. Therefore, $\q_{1,j+1}$ is reachable by the sequence of moves
\begin{gather*}
\q_{1,j-2}=(-2a-2b,-a-b)\xrightarrow{(+a,+b)}(-a-2b,-a)\xrightarrow{(+a,+b)}(-2b,b-a)\\
\xrightarrow{(+a,+b)}(a-2b,2b-a)\xrightarrow{(+b,+a)}(a-b,2b)\xrightarrow{(+b,+a)}(a,a+2b)\xrightarrow{(+b,+a)}\q_{1,j+1}.
\end{gather*}
\item[\textit{Case $2(c)$}:] $2b\mod{0}{n}$ and $\q_{1,j-1}=(0,-a-b)$.\\
Note that $j>1$ since $\q_{1,1-1}=(a,b)\neq(0,-a-b)$. Also note that $-2a-2b\equiv-2a\not\equiv-2b\mod{0}{n}$ and $-2a-b\equiv -2a+b\not\mod{0}{n}$. Therefore, $\q_{1,j+1}$ is reachable by the sequence of moves
\begin{gather*}
\q_{1,j-2}=(-a-b,-2a-2b)\xrightarrow{(+a,+b)}(-b,-2a-b)\xrightarrow{(+a,+b)}(a-b,-2a)\\
\xrightarrow{(+a,+b)}(2a-b,b-2a)\xrightarrow{(+b,+a)}(2a,b-a)\xrightarrow{(+b,+a)}(2a+b,b)\xrightarrow{(+b,+a)}\q_{1,j+1}.
\end{gather*}
\end{itemize}
\end{itemize}

After showing that $\q_{1,j}$ is reachable for all $0\leq j\leq n-1$, we will prove by induction on $i$ that $\q_{i,j}$ is reachable for all $2\leq i\leq n-1$ and $0\leq j\leq n-1$. Assume that for some $1\leq i\leq n-2$, $\q_{i,j}$ is reachable for all $0\leq j\leq n-1$. If $\q_{i,j}$ is not a winning position, then $\q_{i+1,j}$ is reachable by the move
$$\q_{i,j}\xrightarrow{(+a,+b)}\q_{i+1,j}.$$
Otherwise, if $\q_{i,j}$ is a winning position, i.e., $\q_{i,j}=(ia+j(a+b),ib+j(a+b))\in\{(0,i(b-a)),(i(a-b),0)\}$, then we proceed by considering the following cases.
\begin{itemize}
\setlength{\itemindent}{15pt}
\item[\textit{Case $1$}:] $\q_{i,j-1}$ is not a winning position.
\begin{itemize}
\setlength{\itemindent}{30pt}
\item[\textit{Case $1(a)$}:] $\q_{i+1,j-1}$ is not a winning position.\\
By \eqref{eqn:k(a+b)}, $(i+1)(b-a)\not\mod{0}{n}$. Hence, $\q_{i+2,j-1}\in\{(a-b,(i+1)(b-a)),((i+1)(a-b),b-a)\}$ is not a winning position. Therefore, $\q_{i,j+1}$ is reachable by the sequence of moves
$$\q_{i,j-1}\xrightarrow{(+a,+b)}\q_{i+1,j-1}\xrightarrow{(+a,+b)}\q_{i+2,j-1}\xrightarrow{(+b,+a)}\q_{i+1,j}.$$
\item[\textit{Case $1(b)$}:] $\q_{i+1,j-1}$ is a winning position.\\
Since $\q_{i+1,j-1}\in\{(-b,-a+i(b-a)),(-b+i(a-b),-a)\}$, we have $\q_{i+1,j-1}\in\{(-b,0),(0,-a)\}$. Then $\q_{i,j-2}\in\{(-2a-2b,-a-2b),(-2a-b,-2a-2b)\}$ and $\q_{i+2,j-2}\in\{(-2b,-a),(-b,-2a)\}$.

\begin{itemize}
\setlength{\itemindent}{45pt}
\item[\textit{Case $1(b)(i)$}:] $\q_{i,j-2}$ and $\q_{i+2,j-2}$ are not winning positions.\\
Note that $\q_{i+1,j-2}\in\{(-a-2b,-a-b),(-a-b,-2a-b)\}$, $\q_{i+3,j-2}\in\{(a-2b,b-a),(a-b,b-2a)\}$, and $\q_{i+2,j-1}\in\{(a-b,b),(a,b-a)\}$ are not winning positions. Therefore, $\q_{i,j+1}$ is reachable by the sequence of moves
\begin{gather*}
\q_{i,j-2}\xrightarrow{(+a,+b)}\q_{i+1,j-2}\xrightarrow{(+a,+b)}\q_{i+2,j-2}\\
\xrightarrow{(+a,+b)}\q_{i+3,j-2}\xrightarrow{(+b,+a)}\q_{i+2,j-1}\xrightarrow{(+b,+a)}\q_{i+1,j}.
\end{gather*}
\item[\textit{Case $1(b)(ii)$}:] $\q_{i+2,j-2}=(-2b,-a)$ is a winning position, i.e., $2b\mod{0}{n}$ and $\q_{i+2,j-2}=(0,-a)$.\\
Since $b<n$ and $n$ divides $2b$, we have $n=2b$, which is an even number. This implies that $n>3$, thus $-3a-b\equiv-3(a+b)\not\mod{0}{n}$ by \eqref{eqn:k(a+b)}. Moreover, $-2a-b\equiv-2a+b\not\mod{0}{n}$ by the given conditions, and $-2a\not\mod{0}{n}$ since $a<b=\frac{n}{2}$. Therefore, $\q_{i+1,j}$ is reachable by the sequence of moves
\begin{gather*}
\q_{i,j-3}=(-3a-b,-2a-b)\xrightarrow{(+a,+b)}(-2a-b,-2a)\xrightarrow{(+a,+b)}(-a-b,b-2a)\\
\xrightarrow{(+a,+b)}(-b,-2a)\xrightarrow{(+a,+b)}(a-b,b-2a)\xrightarrow{(+b,+a)}(a,b-a)\\
\xrightarrow{(+b,+a)}(a+b,b)\xrightarrow{(+b,+a)}\q_{i+1,j}.
\end{gather*}
\item[\textit{Case $1(b)(iii)$}:] $\q_{i+2,j-2}=(-b,-2a)$ is a winning position, i.e., $2a\mod{0}{n}$ and $\q_{i+2,j-2}=(-b,0)$.\\
Since $a<n$ and $n$ divides $2a$, we have $n=2a$, which is an even number. This implies that $n>3$, thus $-a-3b\equiv-3(a+b)\not\mod{0}{n}$ by \eqref{eqn:k(a+b)}. Moreover, $-a-2b\equiv a-2b\not\mod{0}{n}$ by the given conditions, and $-2b\not\mod{0}{n}$ since $\frac{n}{2}=a<b<n$. Therefore, $\q_{i+1,j}$ is reachable by the sequence of moves
\begin{gather*}
\q_{i,j-3}=(-a-2b,-a-3b)\xrightarrow{(+a,+b)}(-2b,-a-2b)\xrightarrow{(+a,+b)}(a-2b,-a-b)\\
\xrightarrow{(+a,+b)}(-2b,-a)\xrightarrow{(+a,+b)}(a-2b,b-a)\xrightarrow{(+b,+a)}(a-b,b)\\
\xrightarrow{(+b,+a)}(a,a+b)\xrightarrow{(+b,+a)}\q_{i+1,j}.
\end{gather*}
\item[\textit{Case $1(b)(iv)$}:] $\q_{i,j-2}=(-2a-2b,-a-2b)$ is a winning position, i.e., $a+2b\mod{0}{n}$ and $\q_{i,j-2}=(-a,0)$.\\
Note that $n>3$; otherwise, $a=1$ and $b=2$ by $a<b<n$, which contradicts that $a+2b\mod{0}{n}$. By \eqref{eqn:k(a+b)}, $-2a-b\equiv-3(a+b)\not\mod{0}{n}$. Therefore, $\q_{i+1,j}$ is reachable by the sequence of moves
\begin{gather*}
\q_{i,j-3}=(-2a-b,-a-b)\xrightarrow{(+a,+b)}(-a-b,-a)\xrightarrow{(+a,+b)}(-b,b-a)\\
\xrightarrow{(+a,+b)}(a-b,2b-a)\xrightarrow{(+a,+b)}(2a-b,3b-a)=(a-3b,b-2a)\\
\xrightarrow{(+b,+a)}(a-2b,b-a)\xrightarrow{(+b,+a)}(a-b,b)\xrightarrow{(+b,+a)}\q_{i+1,j}.
\end{gather*}
\item[\textit{Case $1(b)(v)$}:] $\q_{i,j-2}=(-2a-b,-2a-2b)$ is a winning position, i.e., $2a+b\mod{0}{n}$ and $\q_{i,j-2}=(0,-b)$.\\
Note that $n>3$; otherwise, $a=1$ and $b=2$ by $a<b<n$, which contradicts that $2a+b\mod{0}{n}$. By \eqref{eqn:k(a+b)}, $-a-2b\equiv-3(a+b)\not\mod{0}{n}$. Therefore, $\q_{i+1,j}$ is reachable by the sequence of moves
\begin{gather*}
\q_{i,j-3}=(-a-b,-a-2b)\xrightarrow{(+a,+b)}(-b,-a-b)\xrightarrow{(+a,+b)}(a-b,-a)\\
\xrightarrow{(+a,+b)}(2a-b,b-a)\xrightarrow{(+a,+b)}(3a-b,2b-a)=(a-2b,b-3a)\\
\xrightarrow{(+b,+a)}(a-b,b-2a)\xrightarrow{(+b,+a)}(a,b-a)\xrightarrow{(+b,+a)}\q_{i+1,j}.
\end{gather*}
\end{itemize}
\end{itemize}

\item[\textit{Case $2$}:] $\q_{i,j-1}=(-a-b,(i-1)b-(i+1)a)$ is a winning position, i.e., $(i-1)b-(i+1)a\mod{0}{n}$ and $\q_{i,j-1}=(-a-b,0)$.\\
If $2b\mod{0}{n}$, then $\q_{i+1,j-1}=(-b,b)=(b,b)\in\{\q_{0,j'}:0\leq j'\leq n-1\}$. This implies that $i=n-1$, violating the bound given in the induction assumption. Hence, $2b\not\mod{0}{n}$.
\begin{itemize}
\setlength{\itemindent}{30pt}
\item[\textit{Case $2(a)$}:] $a+2b\not\mod{0}{n}$.\\
The position $\q_{i+1,j}$ is reachable by the sequence of moves
\begin{gather*}
\q_{i,j-2}=(-2a-2b,-a-b)\xrightarrow{(+a,+b)}(-a-2b,-a)\xrightarrow{(+a,+b)}(-2b,b-a)\\
\xrightarrow{(+a,+b)}(a-2b,2b-a)\xrightarrow{(+b,+a)}(a-b,2b)\xrightarrow{(+b,+a)}\q_{i+1,j}.
\end{gather*}
\item[\textit{Case $2(b)$}:] $a+2b\mod{0}{n}$.\\
Note that $n>3$; otherwise, $a=1$ and $b=2$ by $a<b<n$, which contradicts that $a+2b\mod{0}{n}$. By \eqref{eqn:k(a+b)}, $-2a-b\equiv-3(a+b)\not\mod{0}{n}$. Moreover, $-2a\equiv-a+2b\not\mod{0}{n}$ by the given conditions. Therefore, $\q_{i+1,j}$ is reachable by the sequence of moves
\begin{gather*}
\q_{i,j-3}=(-2a-b,-2a-2b)\xrightarrow{(+a,+b)}(-a-b,-2a-b)\xrightarrow{(+a,+b)}(-b,-2a)\\
\xrightarrow{(+a,+b)}(a-b,b-2a)\xrightarrow{(+a,+b)}(2a-b,2b-2a)\xrightarrow{(+b,+a)}(2a,2b-a)\\
\xrightarrow{(+b,+a)}(2a+b,2b)\xrightarrow{(+b,+a)}\q_{i+1,j}.
\end{gather*}
\end{itemize}

\item[\textit{Case $3$}:] $\q_{i,j-1}=((i-1)a-(i+1)b,-a-b)$ is a winning position, i.e., $(i-1)a-(i+1)b\mod{0}{n}$ and $\q_{i,j-1}=(0,-a-b)$.\\
If $2a\mod{0}{n}$, then $\q_{i+1,j-1}=(a,-a)=(a,a)\in\{\q_{0,j'}:0\leq j'\leq n-1\}$. This implies that $i=n-1$, violating the bound given in the induction assumption. Hence, $2a\not\mod{0}{n}$.
\begin{itemize}
\setlength{\itemindent}{30pt}
\item[\textit{Case $3(a)$}:] $2a+b\not\mod{0}{n}$.\\
The position $\q_{i+1,j}$ is reachable by the sequence of moves
\begin{gather*}
\q_{i,j-2}=(-a-b,-2a-2b)\xrightarrow{(+a,+b)}(-b,-2a-b)\xrightarrow{(+a,+b)}(a-b,-2a)\\
\xrightarrow{(+a,+b)}(2a-b,b-2a)\xrightarrow{(+b,+a)}(2a,b-a)\xrightarrow{(+b,+a)}\q_{i+1,j}.
\end{gather*}
\item[\textit{Case $3(b)$}:] $2a+b\mod{0}{n}$.\\
Note that $n>3$; otherwise, $a=1$ and $b=2$ by $a<b<n$, which contradicts that $2a+b\mod{0}{n}$. By \eqref{eqn:k(a+b)}, $-a-2b\equiv-3(a+b)\not\mod{0}{n}$. Moreover, $-2b\equiv2a-b\not\mod{0}{n}$ by the given conditions. Therefore, $\q_{i+1,j}$ is reachable by the sequence of moves
\begin{gather*}
\q_{i,j-3}=(-2a-2b,-a-2b)\xrightarrow{(+a,+b)}(-a-2b,-a-b)\xrightarrow{(+a,+b)}(-2b,-a)\\
\xrightarrow{(+a,+b)}(a-2b,b-a)\xrightarrow{(+a,+b)}(2a-2b,2b-a)\xrightarrow{(+b,+a)}(2a-b,2b)\\
\xrightarrow{(+b,+a)}(2a,a+2b)\xrightarrow{(+b,+a)}\q_{i+1,j}.
\end{gather*}
\end{itemize}
\end{itemize}
\end{proof}

\section{Proof of Theorem~\ref{thm:theorem2}}\label{sec:section3}

\begin{proof}
Let $d=\gcd(a,b)$ and $\delta=\gcd(m,n)$, and further let $a=da_0$, $b=db_0$, and $n=\delta n_0$ for some integers $a_0$, $b_0$, and $n_0$. Note that $\gcd(d,\delta)=1$ since $\gcd(m,n,a,b)=1$.

Suppose that $m\mid(b^2-a^2)$ and $m\mid\gcd(a,b)\gcd(m,n)$. Then $m=d\delta/c$ for some $c\mid d$, and $(d\delta/c)\mid d^2(b_0^2-a_0^2)$ implies that $\delta\mid cd(b_0^2-a_0^2)$. Since $\gcd(c,\delta)=\gcd(d,\delta)=1$, we have $\delta\mid(b_0^2-a_0^2)$. Let $\delta=\delta^+\delta^-$, where $\delta^+\mid(b_0+a_0)$ and $\delta^-\mid(b_0-a_0)$. Then $b_0=s\delta^-+a_0$ for some integer $s$. Moreover, $\gcd(a_0,\delta^-)=1$ since $\gcd(a_0,b_0)=1$.

We will now show that if $(x_0,0)$ is a reachable winning position, then $x_0\mod{0}{m}$. For any reachable position $(ai+bj,aj+bi)$ with $aj+bi\mod{0}{n}$, we have $da_0j+d(s\delta^-+a_0)i=t\delta n_0$ for some integer $t$. Rearranging the terms, we have $da_0(j+i)=\delta^-(-dsi+t\delta^+n_0)$, so $\delta^-\mid(j+i)$ since $\gcd(da_0,\delta^-)=1$.

As a result, $\delta^+\delta^-\mid(b_0+a_0)(j+i)$, so $\delta\mid(a_0i+b_0j+a_0j+b_0i)$. Recalling that $n\mid(aj+bi)$, we have $\delta\mid(a_0j+b_0i)$. Consequently, $\delta\mid(a_0i+b_0j)$, which implies that $d\delta\mid(ai+bj)$. Therefore, $x_0=ai+bj\mod{0}{m}$, thus proving the sufficient condition for all reachable winning positions being of the form $(0,y)$.

To prove the necessary condition, we assume that all reachable winning positions are of the form $(0,y)$. First, consider the case when $m=a+b$. Then $m\mid(b^2-a^2)$ trivially. Moreover, $d\mid m$ and $\delta\mid m$, which implies that $d\delta\mid m$ since $\gcd(d,\delta)=1$. Hence, $m=\ell d\delta$ for some positive integer $\ell$, or equivalently, $\delta=(a_0+b_0)/\ell$. Assume by way of contradiction that $\ell>1$.

Let $k$ be the smallest positive integer such that $(ka,kb)$ is a reachable winning position. Then $ka=\lcm(a,m)=\lcm(a,a+b)=a(a_0+b_0)$, implying that $k=a_0+b_0$. Thus $\delta<k$, so the positions $(\delta a,\delta b+um)$ are reachable for all $u\geq0$ by the following sequence of moves:
\begin{equation*}
\begin{gathered}\underset{\delta\text{ times of }(+a,+b)}{\underbrace{(0,0)\xrightarrow{(+a,+b)}(a,b)\xrightarrow{(+a,+b)}(2a,2b)\xrightarrow{(+a,+b)}\dotsb\xrightarrow{(+a,+b)}(\delta a,\delta b)}}\\
\hspace{70pt}\left.\begin{array}{c}
\xrightarrow{(+b,+a)}((\delta-1)a,(\delta-1)b+m)\xrightarrow{(+a,+b)}(\delta a,\delta b+m)\\
\vdots\\
\xrightarrow{(+b,+a)}((\delta-1)a,(\delta-1)b+um)\xrightarrow{(+a,+b)}(\delta a,\delta b+um).
\end{array}\right\}\substack{u\text{ times of }(+b,+a)\\\text{ and }(+a,+b)}
\end{gathered}
\end{equation*}
Since $\delta=\gcd(m,n)$, there exist positive integers $u$ and $v$ such that $\delta b=vn-um$. Hence, $(\delta a,\delta b+um)$ is a reachable winning position of the form $(x,0)$ where $x\not\mod{0}{m}$, which is a contradiction. Therefore, $\ell=1$ and $m=\gcd(a,b)\gcd(m,n)$.

It remains to consider the case when $m\neq a+b$. For each positive integer $r$, let $\mathcal{D}_r=\{\p_{r,i}=(a(r-i)+bi,ai+b(r-i)):0\leq i\leq r\}$. Note that $b-a\not\mod{0}{m}$, so for any positive integer $r$ and $0\leq i\leq r$, $a(r-i)+bi$ and $a(r-i-1)+b(i+1)$ are not both congruent to $0$ modulo $m$. In other words, $\p_{r,i}$ and $\p_{r,i+1}$ are not both winning positions. As a result, if both $\p_{r,i}$ and $\p_{r,i+1}$ are reachable positions, then at least one of the moves
$$\p_{r,i}\xrightarrow{(+b,+a)}\p_{r+1,i+1}\text{ and }\p_{r,i+1}\xrightarrow{(+a,+b)}\p_{r+1,i+1}$$
is valid, implying that $\p_{r+1,i+1}$ is reachable.

Note that $\p_{1,0}$, $\p_{1,1}$, $\p_{2,0}$, $\p_{2,1}$, and $\p_{2,2}$ are all reachable. Let $r\geq2$ such that all positions in $\{\p_{r,i}:\sigma\leq i\leq\tau\}$ are reachable for some $0\leq\sigma<\sigma+2\leq\tau\leq r$. Repeatedly applying the previous argument, we see that all positions in 
\begin{equation}\label{reachableset}
\{\p_{r+1,i}:\sigma+1\leq i\leq\tau\}\cup\{\p_{r+2,i}:\sigma+2\leq i\leq\tau\}\cup\{\p_{r+3,i}:\sigma+3\leq i\leq\tau\}
\end{equation}
are reachable. Furthermore, we claim that $\p_{r+2,\sigma+1}$, $\p_{r+2,\tau+1}$, $\p_{r+3,\sigma+1}$, $\p_{r+3,\sigma+2}$, $\p_{r+3,\tau+1}$, and $\p_{r+3,\tau+2}$ are also reachable, and we provide the proof below.

If $\p_{r,\sigma}$ is a winning position, then $\p_{r,\sigma}=(0,y_0)$ for some integer $y_0$. Hence, $\p_{r,\sigma+1}=(b-a,y_0+a-b)$, $\p_{r+1,\sigma+1}=(b,y_0+a)$, and  $\p_{r+2,\sigma+1}=(a+b,y_0+a+b)$ are all reachable non-winning positions, which further implies that both $\p_{r+3,\sigma+1}$ and $\p_{r+3,\sigma+2}$ are reachable.

On the other hand, if $\p_{r,\sigma}$ is not a winning position, then $\p_{r+1,\sigma}$ is reachable. Now, if $\p_{r+1,\sigma}$ is a winning position, then $\p_{r+1,\sigma}=(0,y_1)$ for some integer $y_1$. Hence, both $\p_{r+1,\sigma+1}=(b-a,y_1+a-b)$ and $\p_{r+2,\sigma+1}=(b,y_1+a)$ are reachable non-winning positions, thus both $\p_{r+3,\sigma+1}$ and $\p_{r+3,\sigma+2}$ are also reachable. Otherwise, if $\p_{r+1,\sigma}$ is not a winning position, then $\p_{r+2,\sigma}$ and $\p_{r+2,\sigma+1}$ are reachable. Recalling from \eqref{reachableset} that $\p_{r+2,\sigma+2}$ is also reachable, it follows that both $\p_{r+3,\sigma+1}$ and $\p_{r+3,\sigma+2}$ are also reachable. Similar arguments will show that $\p_{r+2,\tau+1}$, $\p_{r+3,\tau+1}$, and $\p_{r+3,\tau+2}$ are all reachable, thus concluding our proof for the claim.

Since $\p_{r,i}=\p_{r,i'}$ if $i'=i+\lcm(m,n)$, the positions in $\mathcal{D}_r$ are periodic, meaning that as long as $\tau-\sigma\geq\lcm(m,n)$, we have $\{\p_{r,i}:\sigma\leq i\leq\tau\}=\mathcal{D}_r$. From the claim above, we observe that if all positions in $\{\p_{r,i}:\sigma\leq i\leq\tau\}$ are reachable for some $0\leq\sigma<\sigma+2\leq\tau\leq r$, then all positions in
$$\{\p_{r+1,i}:\sigma+1\leq i\leq\tau\}\cup\{\p_{r+2,i}:\sigma+1\leq i\leq\tau+1\}\cup\{\p_{r+3,i}:\sigma+1\leq i\leq\tau+2\}$$
are also reachable. Applying the claim repeatedly, we know that all positions in
\begin{gather*}
\{\p_{r+3w+1,i}:\sigma+w+1\leq i\leq\tau+2w\}\cup\{\p_{r+3w+2,i}:\sigma+w+1\leq i\leq\tau+2w+1\}\\
\cup\{\p_{r+3w+3,i}:\sigma+w+1\leq i\leq\tau+2w+2\}
\end{gather*}
are reachable for all positive integers $w$. Hence, for all $w>\lcm(m,n)$, all positions in $\mathcal{D}_{r+3w+1}\cup\mathcal{D}_{r+3w+2}\cup\mathcal{D}_{r+3w+3}$ are reachable. Moreover, since $\mathcal{D}_r=\mathcal{D}_{r'}$ if $r'=r+\lcm(m,n)$, we conclude that every position $(ai+bj,aj+bi)$ is reachable.

From this, we see that if $i=\lcm(m,n)-a$ and $j=b$, then $(ai+bj,aj+bi)=(b^2-a^2,0)$ is reachable. Based on the assumption that all winning positions are of the form $(0,y)$, we have $m\mid(b^2-a^2)$. Similarly, letting $i=n$ and $j=\lcm(m,n)$, we know that both $(ai+bj,aj+bi)=(an,0)$ and $(aj+bi,ai+bj)=(bn,0)$ are reachable. Again, since all winning positions are of the form $(0,y)$, we have $m\mid an$ and $m\mid bn$. This implies that $m\mid\gcd(a,b)n$, thus $m\mid\gcd(a,b)\gcd(m,n)$.
\end{proof}

\section{Acknowledgments}
These results are based on work supported by the National Science Foundation under grant numbered DMS-1852378.


\begin{thebibliography}{1}





\bibitem{hkksw}
J.~Harrington, K.~Karhadkar, M.~Kohutka, T.~Stevens, and T.~W.~H.~Wong, Two dependent probabilistic chip-collecting games, \textit{Discrete Appl.\ Math.\/} \textbf{288} (2021), 74--86.

\bibitem{lt}
H.H.~Leung and T.~Thanatipanonda, A probabilistic two-pile game, \textit{J.\ Integer Seq.\/} \textbf{22} (2019), 19.4.8.

\bibitem{lt2020}
H.H.~Leung and T.~Thanatipanonda, Game of pure chance with restricted boundary, \textit{Discrete Appl.\ Math.\/} \textbf{283} (2020), 613--625.





\bibitem{wx}
T.W.H.~Wong and J.~Xu, A probabilistic take-away game, \textit{J.\ Integer Seq.\/} \textbf{21} (2018), 18.6.3.

\end{thebibliography}
\end{document}